\documentclass[12pt]{article}

\usepackage{graphicx, float, ulem, subcaption, amssymb, lipsum, graphics, mathtools, amsmath, amsfonts}
\usepackage[colorlinks,linkcolor=blue,citecolor=blue,urlcolor=blue]{hyperref}
\usepackage[linesnumbered,ruled,vlined]{algorithm2e}
\usepackage[english]{babel}
\usepackage{epstopdf, xcolor, soul, framed}
\usepackage[numbers]{natbib}

\SetCommentSty{mycommfont}
\SetKwInput{KwInput}{Input} 
\SetKwInput{KwOutput}{Output} 

\graphicspath{{images/}}

\colorlet{shadecolor}{yellow}

\sethlcolor{yellow}

\begin{document}


\title{Understanding Work/Life Balance through the Mathematical Model}

\author{
  Francy Singla\footnotemark[2] \hspace{0.1 cm}, 
  \ Vivek Kumar\footnotemark[1] \\
  \vspace{0.25cm} 
  \normalsize \textit{Department of Applied Mathematics, Delhi Technological University, Delhi, India-110042}
}

\date{} 
\maketitle
\renewcommand{\thefootnote}{\fnsymbol{footnote}}

\footnotetext[2]{E-mail: francysingla1606@gmail.com}
\footnotetext[1]{Corresponding author. E-mail: vivekkumar.ag@gmail.com}

{\let\thefootnote\relax\footnotetext{\textbf{Note to Editorial Team:} The LaTeX source files and image directory required for typesetting can be downloaded from this secure repository: \url{https://drive.google.com/file/d/1Y2Xoh4hW_vuIEaF8ljfrt3MdYDfct-Tc/view?usp=sharing}}
}


\begin{abstract}
Work/life balance (WLB) remains a central yet elusive concept in social and organizational sciences, often described qualitatively without a unified quantitative framework. To address this gap, a novel compartmental mathematical model is proposed that formalizes WLB as a dynamic population-level phenomenon. The individuals are categorized into 4 compartments, namely, regularly employed, temporarily employed, imbalanced (infected), and balanced (recovered), representing balanced and imbalanced states. Role-based complaints govern the transitions between states, the conflict conversion rates of these complaints, and the recovery rate. The model, structured as a system of nonlinear ordinary differential equations, is analyzed for the existence of the equilibrium and its stability. In the analysis, the stable disease-free equilibrium is found to exist, representing complete societal balance. The other endemic equilibrium is shown to be infeasible through simulations and thus requires certain parameter constraints, which represent a specific scenario that may occur in the population. However, it is found that the endemic equilibrium always remains unstable. The model's ability to quantify how micro-level stressors spread into macro-level imbalance is demonstrated by numerical simulations, which also validate the analytical results. This work establishes a rigorous mathematical foundation for WLB, enabling future empirical and mathematical calibrations.
\end{abstract}



\section{Introduction}
\label{Intro}

Work-life balance has been a subject of research and concern among many researchers in the fields of human resource management (HRM), psychology, social sciences, etc. The topic remained in the shadows for many years and has started earning the spotlight over the past decade. The HRM started recognizing the requirement and demand of their employees for the work/life balance, with which the need for the development of mechanisms that facilitate WLB came forth. Although numerous scholars have proposed various definitions of WLB over the past few decades, no definition has been agreed upon by WLB scholars as the gold standard. The meaning of WLB is clear in the name itself, that is, maintaining equilibrium between one's professional and personal lives, or more specifically, between one's work and non-work life. It involves meeting work obligations while nurturing personal life. Prioritizing one segment of life while ignoring the other causes an imbalance between the two, which has severe consequences. 

The absence of balance, i.e., the imbalance, occurred as a consequence of a shift in the work culture over the past years, along with the change in lifestyle of the people. The present work culture demands employees to work for longer durations for good pay and promotes the workaholic behavior to achieve their desired lifestyle, goals, success, and fame. This has initiated an unhealthy cycle in people's lives where they either work endlessly, ignoring their health, family, personal life, and other aspects of life. The brain remains occupied almost all the time, keeping it under stress and awarding little to no time to relax. It leads to people feeling overwhelmed and experiencing cognitive overload when presented with more information and tasks than their brains can process. Due to cognitive overload, the brain struggles to make space for any other task, emotion, or mindful presence, which has a huge behavioral and psychological impact. This creates a pseudo-sense in people's minds of a lack of time in their lives. As a result, they assume that the work/life balance is just a measure of time and that dividing equal time between both work and non-work life balances both.

While the spotlight usually remains on the effect of poor WLB on one's non-work life due to overindulgence in work life, it becomes important to understand that the vice versa also exists. This means that the poor WLB also affects one's work life due to overindulgence in non-work life. Thus, the concept of WLB is bidirectional   (\citet{frone1992antecedents}). However, the question remains the same: what defines balance in one's work and non-work life? How does one know if the WLB exists in one's life—are they in a state of balance or imbalance? Is the state of balance continuous? 

Before proceeding to understand WLB, it is also important to know what it does not account for. Firstly, the balance between work and non-work life does not equate to equality; that is, trying to allocate the same amount of time to all of your job and personal obligations is unreasonable   (\citet{meenakshi2013importance}). Secondly, one's ideal work/life balance changes over time, frequently every day, meaning that the ideal balance today will most likely change from the ideal balance tomorrow   (\citet{meenakshi2013importance}). Further, there is no ideal equilibrium that applies to everyone, because all of us have various objectives and lives; hence, the ideal work-life balance varies from person to person   (\citet{meenakshi2013importance}). This again raises the question: Is there any way to measure collective work-life balance for all? 

In the existing literature, researchers in the past attempted to develop multiple scales to measure different components that account for WLB. In their paper, \citet{article} discuss three studies they conducted to analyze different measures proposed by different researchers and further gave a scale of their own rooted in three forms of role conflict and dual direction of WLB. The Role Theory conceptualizes WLB as a state of minimal inter-role conflict accompanied by satisfactory functioning across work and non-work domains. However, this understanding is still primarily descriptive, which makes it difficult to test, anticipate, or systematically intervene at scale. \citet{gulmez2024heuristic} develop a mathematical model to prioritize employees’ self-selection for their shifts and further increase the self-selection rates of working time and starting hours with a special interest to the work/life balance of healthcare workers. They discuss that a clear gap exists in modeling the problem mathematically and how the definition, modeling, and optimization of this problem have been a missing part of the literature. The problem should be well-defined in mathematical structures so as to make it ready to solve.  

The absence of a quantitative framework limits our ability to assess how individual experiences of role strain contribute to population-level trends of stress, declining mental health, and social unrest. This paper addresses this gap by introducing a mathematical model that treats WLB as a dynamic, population-level process. The model acts as a tool to understand the WLB of a single individual and a group of individuals (population). How long an individual remains in the state of imbalance, the factors that contribute to an increase in balance or imbalance, how long it takes to return to the state of balance, etc., all these questions become answerable with the mathematical model. Thus, WLB becomes a predictive, analyzable construct when viewed through the prism of dynamical systems, moving beyond static definitions.

\textit{Motivation:} The inspiration of the study came from the lack of answers to the question: Is WLB just a function of time spent at work? Or how can one know their work life is in balance with their non-work life, or is it imbalanced? Is there a mathematical formula (quantitative tool) to measure the balance? Moreover, as work/life imbalance intensifies, individual stress levels rise, and the ability to empathize with people around them decreases. Additionally, it can cause feelings of loneliness, despondency, and separation, which can further lead to anxiety and depression. This results in a crisis of mental health in the entire population, telling us why there is an urgent need to pay attention to this new widespread problem before it becomes an infectious yet non-communicable epidemic that affects the entire human race.

The paper is structured as follows:  Section \ref{Pre} discusses the definitions, concepts, and prerequisite information to understand the model. Section \ref{Model} gives the detailed description of the formulation of the dynamical WLB model. In section \ref{EA}, we conduct the equilibrium analysis, while section \ref{SA} investigates the local stability of these equilibria. Section \ref{numerical sim} validates the theoretical results through numerical simulations under varied parameter regimes. Finally, section \ref{Con} discusses the results and interpretations made from the study of the model and discusses their implications for theory and practice.

\section{Preliminaries}
\label{Pre}

The balance between work and family was the main focus of research in previous decades, as in those days, non-work life was equivalent to family life. In contrast, today's non-work life includes more than just family life; it also includes friends, hobbies, self-care, and other pursuits. Thus, scholars such as \textit{Greenhaus, Beutall, and Powell} described the notions using the terms 'work/family conflict' (WFC) instead of 'work/life conflict' and 'work/family enrichment' (WFE) instead of 'work/life enrichment'   (\citet{greenhaus1985sources}, \citet{greenhaus2006work}). Researchers in the 2000s started adopting more general terminology, substituting "life" for "family" and defining concepts in a broader sense, recognizing that family is only a part of life outside of work.

The foundation of the literature of WLB is based on the \textit{role theory}. According to the role theory, a person's life is made up of various roles, such as those they play at work (like supervisor) and outside of it (like spouse, friend).  Role conflict arises when two or more sets of pressures are present at the same time, and meeting one set makes meeting the other set more challenging   (\citet{fisher2003modelling}). Furthermore, the dominance of the conflict perspective in the work–family literature is rooted in scarcity theory, which assumes that the personal resources of time, energy, and attention are finite, and that the devotion of greater resources to one role necessitates the devotion of lesser resources to the other role   (\citet{GREENHAUS2003291}). Thus, individuals who participate in both work and family roles are likely to experience conflict between these roles. Combining these theories would imply that juggling several responsibilities creates competition for the same resources, which might result in stressors like health issues or job discontent   (\citet{fisher2003modelling}). 

\textbf{Work/Family Conflict:} A form of inter-role conflict in which the role pressures from work and family domains are mutually incompatible in some respect   (\citet{greenhaus1985sources}).
In the broader sense, \textit{\textbf{work/life conflict}} is defined as a type of inter-role conflict where the demands of work and personal life are incompatible.
The types of conflicts that may occur in multiple roles are further smartly divided into three categories by \textit{Greenhaus and Beutall}, which are namely time-based conflict, strain-based conflict, and the behavioral-based conflict   (\citet{mcmillan2011constructs}).

These constructs of conflict, balance, and enrichment are the cross-domain effects, meaning that the experiences and decisions in one domain are capable of influencing outcomes in the other domain   (\citet{mcmillan2011constructs}). This implies that the concepts of WLB, WFC, and WFE are bidirectional; that is, conflicts, balance, and enrichment in roles from work life can interfere with the roles outside the work and vice versa. In their review paper \citet{sirgy2018work} discuss multiple definitions of WLB proposed by the researchers, which they categorized into two key dimensions: 1. Role engagement in multiple roles in work and non-work life, and 2. Minimal conflict between work and non-work roles. The key characteristics of role engagement are equal time and involvement and balanced satisfaction, while those of minimal conflict are role enrichment with no role conflict and management of resources to minimize role conflicts.

In our study, we focus on incorporating the second dimension in our model, that is, the minimal role conflict, as the engagement in multiple roles remains an abstract concept, which we leave as a future scope of the study for now. In the following section, we proceed to construct the mathematical model for WLB using the system of ordinary differential equations.

\section{The Mathematical Model for WLB}
\label{Model}

\subsection{Description of the WLB Model:}

The population susceptible to experiencing work/life imbalance is the working population, which is further classified into two categories, namely, the \textit{regularly employed $(S_r)$} and the \textit{temporarily employed $(S_T)$}. The mathematical model given by \citet{misra2011mathematical} discusses how the population moves between these two classes. We use the subscripts 'r' and 'T' to differentiate between regularly and temporarily employed populations, respectively. The division of population recognizes that temporarily employed and regularly employed members of the population experience different amounts of conflicts and balance, as the temporarily employed people are employed for a defined period of time, in contrast to the regularly employed, who are assumed to be employed for an infinite period of time. \newline

Once any member of the working population (both categories) experiences an imbalance between their work and non-work life, that member moves to the \textit{Imbalanced/Infected (I)} compartment, and further, when the members from the imbalanced category restore the balance between work and non-work life, they move to the \textit{Balance Recovered (R)} compartment. The course of movements of the population from one compartment to the other is described below: \newline

Based on the theory discussed in section \ref{Pre}, we take the number of roles and conflicts as the parameters in our model. Let $n$ be the average number of roles every member of the working population plays in their work and non-work life. Let the average number of complaints received per role per member of the population per unit time be given by $c_r$ and $c_T$ for regularly and temporarily employed, respectively. Therefore, $n c_T \boldsymbol{S_T}$ and $n c_r \boldsymbol{S_r}$ give the total number of complaints in all the roles per unit time for temporarily and regularly employed people, respectively.  \newline

The purpose of introducing the complaint parameters ($c_r$ and $c_T$) is to quantify the theory of conflicts. The starting point of a conflict in any role in work or non-work life is the complaints, which are an expression of dissatisfaction. When these complaints in any role remain unaddressed and unresolved, dissatisfaction increases, leading to conflicts. Further, these conflicts fall into any one of the three categories discussed in section \ref{Pre}. In our study, for the simplification of the model, we do not dive deeper into the categories and stick to a broader term that covers all types of conflicts. \newline 

Assume that only a proportion of the total complaints that a member of the population receives in all the roles per unit time ($n c_T \boldsymbol{S_T}$ and $n c_r \boldsymbol{S_r}$) turn into conflicts. Symbolically, we denote them by $q_r$ and $q_T$ $(0 \leq q \leq 1)$, where the subscripts have their usual meanings. Thus, $q_T n c_T \boldsymbol{S_T}$  gives the total number of conflicts in all the roles per unit time among the temporarily employed population, and $q_r n c_r \boldsymbol{S_r}$ gives the total number of conflicts in all the roles per unit time among the regularly employed population. \newline 

The working population leaves their respective susceptible compartment and enters the imbalanced/infected compartment due to these factors. The imbalanced population adopts measures to restore balance between their work and non-work life, say the recovery rate is $\gamma$. Thus, the population leaves the imbalanced/infected compartment and enters the balance recovered compartment by the factor $\gamma \boldsymbol{I}$. \newline

The parameter $S_{ra}$ represents the total number of vacancies for regular employment and $k_1$ the constant of proportionality, both of which are taken from the unemployment model given in \citet{misra2011mathematical}. The movement between the regularly and temporarily employed compartments is given by the factor $k_1(S_{ra} - S_r)$. Refer to Table \ref{table2} for list all the parameters, their symbols, and definitions used in the WLB model. \newline

\textit{\textbf{\uline{ Assumptions:}}} \newline
1. The total population remains constant (N) at any time 't' throughout the study; that is, $N(t)$ = $S_r (t)$ + $S_T (t)$ + $I(t)$ + $R(t)$ \newline
2. The population in $S_r$ and $S_T$ compartments are initially in the state of balance. \newline
3. Once any member reaches the balance-recovered compartment, it is assumed that they remain in the state of balance forever. \newline
4. The average number of roles for both the temporarily and regularly employed classes is the same. \newline
5. All the chosen parameters are real-valued and non-negative, that is, $n,c_r, c_T, q_r, q_T, \gamma, k_1, S_{ra} \geq 0$. \newline
6. On average, every individual indulges in atleast 2 roles, one in work life and the other in non-work life. 

\subsection{The deterministic WLB Model}

\begin{align}
\frac{d\boldsymbol{S_T}}{dt} &= - k_1 (S_{ra} - \boldsymbol{S_r})\boldsymbol{S_T} - q_T c_T n \boldsymbol{S_T} \\
\frac{d\boldsymbol{S_r}}{dt} &=  k_1 (S_{ra} - \boldsymbol{S_r})\boldsymbol{S_T} - q_r c_r n \boldsymbol{S_r} \\
\frac{d\boldsymbol{I}}{dt} &= q_T c_T n \boldsymbol{S_T}+ q_r c_r n \boldsymbol{S_r} - \gamma \boldsymbol{I}  \\ 
\frac{d\boldsymbol{R}}{dt} &= \gamma \boldsymbol{I} 
\end{align}

We reduce this system of 4 differential equations to a system of 3 differential equations by substituting the value of $S_T$ in the above. Since we are studying an isolated system, the total population remains conserved at each point in time; that is, $S_r(t) + S_T(t) + I(t) + R(t) = N(t)$. Substitute the value of $S_T$ as $S_T(t) = N(t) - S_r(t) - I(t) - R(t)$ in eqs. (2) and (3). This transforms our system of 4 differential equations into a system of 3 differential equations. The derivation of the 3-ordered system is as follows:

Eq. (2) transforms to:
\begin{align}
    \frac{d\boldsymbol{S_r}}{dt} &= k_1(S_{ra} - \boldsymbol{S_r})(N - \boldsymbol{S_r} - \boldsymbol{I} - \boldsymbol{R}) - q_rc_rn\boldsymbol{S_r} \\
                    &= k_1(S_{ra} - \boldsymbol{S_r})(N - \boldsymbol{I} - \boldsymbol{R}) - k_1S_{ra}\boldsymbol{S_r} + k_1\boldsymbol{S_r}^2 - q_rc_rn\boldsymbol{S_r} \\
                    &= k_1\boldsymbol{S_r}^2 + k_1S_{ra}(N - \boldsymbol{I} - \boldsymbol{R}) - k_1\boldsymbol{S_r}(N - \boldsymbol{I} - \boldsymbol{R}) - k_1S_{ra}\boldsymbol{S_r} - q_rc_rn\boldsymbol{S_r} \\
                    &= k_1\boldsymbol{S_r}^2 - \boldsymbol{S_r} [k_1(N-\boldsymbol{I} - \boldsymbol{R}) + k_1 S_{ra} + q_rc_rn] + k_1 S_{ra}(N-\boldsymbol{I} - \boldsymbol{R})
\end{align}

and eq. (3) becomes:

\begin{align}
    \frac{d\boldsymbol{I}}{dt} &= q_Tc_Tn(N - \boldsymbol{S_r} - \boldsymbol{I} - \boldsymbol{R}) + q_rc_rn\boldsymbol{S_r} - \gamma\boldsymbol{I} \\
                    &= (q_r c_r n - q_T c_T n) \boldsymbol{S_r} + q_T c_T n (N - \boldsymbol{R}) - (\gamma + q_T c_T n)\boldsymbol{I}
\end{align}

Thus, the reduced WLB system is given by eqs. (11), (12), and (13) below:

\begin{align}
\dfrac{d\boldsymbol{S_r}}{dt} &=  k_1\boldsymbol{S_r}^2 - \boldsymbol{S_r} [k_1(N-\boldsymbol{I} - \boldsymbol{R}) + k_1 S_{ra} + q_rc_rn] + k_1 S_{ra}(N-\boldsymbol{I} - \boldsymbol{R}) \\
\dfrac{d\boldsymbol{I}}{dt} &= (q_r c_r n - q_T c_T n) \boldsymbol{S_r} + q_T c_T n (N - \boldsymbol{R}) - (\gamma + q_T c_T n)\boldsymbol{I}  \\ 
\dfrac{d\boldsymbol{R}}{dt} &= \gamma \boldsymbol{I} 
\end{align}

\subsection{Positive Invariance of the feasible region under the flow of the system}

\textbf{Lemma: }The set $\Omega = {(S_T,S_r,I,R) \in \boldsymbol{R^4}_{\geq0}: S_T + S_r + I + R = N}$ is positively invariant under the flow of the WLB system given by eqs. (1) - (4). \newline

\textbf{\textit{Proof:}} Let $y(t)$ be a vector whose components are the WLB system variables such that $y(t) = [S_T(t), S_r(t), I(t), R(t)]^T$. 
The WLB system can be rewritten as  $y'(t) = f(y(t))$
where $f(y(t)) = \left(\begin{array}{c} 
      - k_1 (S_{ra} - \boldsymbol{S_r})\boldsymbol{S_T} - q_T c_T n \boldsymbol{S_T} \\
     
      k_1 (S_{ra} - \boldsymbol{S_r})\boldsymbol{S_T} - q_r c_r n \boldsymbol{S_r} \\
      
     q_T c_T n \boldsymbol{S_T}+ q_r c_r n \boldsymbol{S_r} - \gamma \boldsymbol{I}  \\

     \gamma \boldsymbol{I}
\end{array} \right)$. 

Obviously, the feasible region $\Omega$ is a closed set, and we want to prove that if the initial condition $f(0) \in \Omega$, then the solution $f(t) \in \Omega$ for all t $\geq 0$. 

The boundary $\partial\Omega$ consists of five hyperspaces where the variables reach their physical limits, which are given below:
\begin{enumerate}
    \item $P_1: R = 0$ 
    \item $P_2: I = 0$
    \item $P_3: S_r = 0$
    \item $P_4: S_T = 0 $
    \item $P_5: S_T + S_r + I + R = N$
\end{enumerate}

Corresponding to the above hyperspaces $P_1, P_2, P_3, P_4, P_5$, take the normal outward vectors (which represent the boundary lines) $w_1 = (0,0,0,-1), w_2 = (0,0,-1,0), w_3 = (0,-1,0,0), w_4 = (-1,0,0,0), w_5 = (1,1, 1, 1)$  respectively. \newline
We show that no solution exits the region $\Omega$, that is, at every boundary, the flow of the system is directed inward or is tangent to the boundary. \newline 
This is shown by calculating the dot product of the normal vectors with $f(y(t))$. If the dot product $f(y(t)) . w_i \leq 0$ for all $i$, then $f(y(t)) \in \Omega$ for all $t\geq 0$. 

For Boundary $P_4$, the dot product
\begin{center}
$w_4 . f(y(t)_{S_T = 0} = (-1,0,0,0) . (\frac{d\boldsymbol{S_T}}{dt}, \frac{d\boldsymbol{S_r}}{dt}, \frac{d\boldsymbol{I}}{dt}, \frac{d\boldsymbol{R}}{dt}) = - \frac{d\boldsymbol{S_T}}{dt}$ 
\end{center}
Substituting the value in the above eq. we get, 
\begin{center}
    $w_4 . f(y(t)_{S_T = 0} =  k_1 (S_{ra} - \boldsymbol{S_r})\boldsymbol{S_T} + q_T c_T n \boldsymbol{S_T} = 0 $
\end{center}
This implies that the flow is tangent to the boundary and $S_T$ cannot become negative.

Similarly, for Boundary $P_3, P_2$ and $P_1$, the flow points into the region or stays on the boundary:
\begin{align}
  w_3 . f(y(t)_{S_r = 0} &= - k_1 (S_{ra} - \boldsymbol{S_r})\boldsymbol{S_T} - q_r c_r n \boldsymbol{S_r} = - k_1 S_{ra}\boldsymbol{S_T} \leq 0 \\
  w_2 . f(y(t)_{I = 0} &=  -q_T c_T n \boldsymbol{S_T} - q_r c_r n \boldsymbol{S_r} + \gamma \boldsymbol{I} = -q_T c_T n \boldsymbol{S_T} - q_r c_r n \boldsymbol{S_r} \leq 0 \\
  w_1 . f(y(t)_{R = 0} &= - \gamma \boldsymbol{I} \leq 0
\end{align}

At last, we show that the total population remains constant in $P_5$.
\begin{center}
    $w_5 . f(y(t))_{S_T + S_r + I + R = N} = \frac{d\boldsymbol{S_T}}{dt} + \frac{d\boldsymbol{S_r}}{dt} + \frac{d\boldsymbol{I}}{dt} + \frac{d\boldsymbol{R}}{dt} = 0$
\end{center}
    This completes the proof that the set $\Omega$ remains invariant under the flow of the WLB system.
\begin{table}[h!]
\centering
\caption{List of parameters used in the WLB Model}
\resizebox{15 cm}{!}{
\begin{tabular}{|c|c|c|c|c|c|c|} 
\hline
\textbf{S. No.} & \textbf{Parameter} & \textbf{Definition}  \\ \hline
1. & n        & average number of roles per member of population      \\ \hline
2. & $c_r$    & average number of complaints per role per regularly employed member per unit time        \\ \hline
3. & $c_T$    & average number of complaints per role per temporarily employed member per unit time   \\ \hline
4. & $q_r$    & conversion rate of complaints to conflicts in regularly employed ($q_r \in [0,1]$)     \\ \hline
5. & $q_T$    & conversion rate of complaints to conflicts in temporarily employed ($q_T \in [0,1]$)               \\ \hline
6. & $\gamma$ & balance recovery rate       \\ \hline
7. & $k_1$    & proportionality constant        \\ \hline
8. & $S_{ra}$ & total number of vacancies for regular employment, ($S_{ra} \geq \boldsymbol{S_r}$)     \\ \hline 
\end{tabular}}
\label{table2}
\end{table}

\section{Equilibrium Analysis}
\label{EA}
Now we perform the equilibrium analysis on the modified system given by equations (11), (12), and (13). For background on non-linear dynamical systems and stability analysis, see \citet{vaidya2026analyzing} and \citet{kulkarni2025numerical}.

The system has two non-negative equilibrium points, which are obtained by solving the algebraic equations:

\begin{align}
k_1\boldsymbol{S_r}^2 - \boldsymbol{S_r} [k_1(N-\boldsymbol{I} - \boldsymbol{R}) + k_1 S_{ra} + q_r c_r n] + k_1 S_{ra}(N-\boldsymbol{I} - \boldsymbol{R}) &= 0\\
(q_r c_r n - q_T c_T n) \boldsymbol{S_r} + q_T c_T n (N - \boldsymbol{R}) - (\gamma + q_T c_T n)\boldsymbol{I} &=0 \\ 
\gamma \boldsymbol{I}  &= 0
\end{align}

From equation (19), we get 
\begin{equation}
    \boldsymbol{I} = 0
\end{equation} 
which we substitute in equations (17) and (18).  \newline
From Eq. (18), we get
\begin{equation}
    \boldsymbol{R} = \frac{\boldsymbol{S_r}(q_r c_r - q_T c_T)}{q_T c_T} + N.
\end{equation}

Substituting the values of $\boldsymbol{I}$ and $\boldsymbol{R}$ from Eq. (20) and (21) into Eq. (17), we get a quadratic equation in $\boldsymbol{S_r}$

\begin{equation}
    (\frac{k_1 q_r c_r}{q_T c_T})\boldsymbol{S_r}^2 - (q_r c_r n + \frac{k_1 S_{ra} q_r c_r}{q_T c_T})\boldsymbol{S_r} = 0
\end{equation}

By solving the above quadratic equation, we get the following.

\begin{align}
\boldsymbol{S_r} &= 0 \\
\boldsymbol{S_r} &= \frac{k_1 S_{ra} + q_T c_T n}{k_1}
\end{align}

Substituting values from Eqs. (23) and (24) in Eq. (18), respectively, we get

\begin{align}
    \boldsymbol{R} &= N \\
    \boldsymbol{R} &= \frac{k_1 S_{ra} q_r c_r - k_1 S_{ra} q_T c_T + n q_T c_T q_r c_r - n (q_T c_T)^2 + k_1 q_T c_T N}{k_1 q_T c_T}
\end{align}

Hence, the two equilibrium points that we have obtained are as follows:
\begin{align}
    E_1 &= (0,0, N) \\
    E_2 &= (\frac{k_1 S_{ra} + q_T c_T n}{k_1}, 0, \frac{k_1 S_{ra} q_r c_r - k_1 S_{ra} q_T c_T + n q_T c_T q_r c_r - n (q_T c_T)^2 + k_1 q_T c_T N}{k_1 q_T c_T})
\end{align}
By substituting the values of $S^*_r$, $I^*$, $R^*$ into the expression $N(t) = S^*_r + S^*_T + I^* + R^*$, the value of $S^*_T$ can be obtained for both $E_1$ and $E_2$, which will give the equilibrium points corresponding to the system given by Eqs. (1)-(4). It is obvious that for $E_1$, the value of $S^*_T$ is 0.

\section{Stability Analysis}
\label{SA}
We investigate the local stability of the model at the two equilibrium points $E_1$ and $E_2$ given by Eq. (18) and (19). The Jacobian matrix for model (5)-(7) is given by: 
\[
\begin{pmatrix}
  2 S_r k_1 - S_{ra} k_1 - c_r n q_r - k_1 (N - I - R) & S_r k_1 - S_{ra} k_1 & S_r k_1 - S_{ra} k_1 \\
  c_r n q_r - c_T n q_T & -c_T n q_T - \gamma & -c_T n q_T \\
  0 & \gamma & 0   
\end{pmatrix}
\] 
The value of the Jacobian at the equilibrium point $E_1$, namely $J_1$, is:
\[
\begin{pmatrix}
  - S_{ra} k_1 - c_r n q_r & - S_{ra} k_1 & - S_{ra} k_1 \\
  c_r n q_r - c_T n q_T & -c_T n q_T - \gamma & -c_T n q_T \\
  0 & \gamma & 0   
\end{pmatrix}
\]    
The eigenvalues of the matrix $J_1$ are $- \gamma, - c_Tnq_T, - c_rnq_r$. All of these eigenvalues are real-valued, as each parameter value is real, and all of them are negative. This implies that the equilibrium point $E_1 = (0, 0, N)$ is asymptotically stable and is a stable node.

The value of the Jacobian at the equilibrium point $E_2$, namely $J_2$, is:
\[
\begin{pmatrix}
c_T n q_T + \dfrac{c_r q_r k_1}{q_T c_T} & c_T n q_T & c_T n q_T \\
c_r n q_r - c_T n q_T & -c_T n q_T - \gamma & -c_T n q_T \\
0 & \gamma & 0   
\end{pmatrix}
\]
The eigenvalues of the matrix $J_2$ are $\dfrac{S_{ra} k_1 c_r q_r \pm \sqrt{S^2_{ra}c^2_r k^2_1 q^2_r + 4 S_{ra} c_r c^2_T k_1 n q^2_T q_r + 4 c_r c^3_T n^2 q^3_T q_r}}{2 c_T q_T}$, and $- \gamma$.
\newline Since every parameter taken in the model is real valued and non-negative, the eigenvalue $\dfrac{ S_{ra} k_1 c_r q_r + \sqrt{S^2_{ra}c^2_r k^2_1 q^2_r + 4 S_{ra} c_r c^2_T k_1 n q^2_T q_r + 4 c_r c^3_T n^2 q^3_T q_r}}{2 c_T q_T}$ remains always positive, making $E_2$ an unstable equilibrium point.

\section{Numerical Simulation}
\label{numerical sim}

In the absence of the actual data, the WLB model (1)-(4) is numerically simulated using MATLAB software. Table \ref{table1} contains the set of parameter values used for the numerical simulations. The parameter values are chosen in such a way that they represent multiple scenarios that may arise in the population. The results from the simulations (1)-(5) are shown through time evolution graphs in figures (1)-(3). A computer simulation is performed for different initial conditions, and the results for that are displayed in figures (4)-(5). The eigenvalues computed for $E_1$ are given in Table \ref{table3}. Further, the phase portraits of the dependent variables are drawn in the figures (6)-(8)

\begin{table}[h!]
\centering
\caption{Parameter values used for numerical simulations}
\resizebox{14 cm}{!}{
\begin{tabular}{|c|c|c|c|c|c|c|c|c|} 
\hline
\textbf{S. No.} & \textbf{Parameter} & \textbf{Simulation 1} & \textbf{Simulation 2} & \textbf{Simulation 3} & \textbf{Simulation 4} & \textbf{Simulation 5} \\ \hline
1. & n        & 2        & 6       &  6      & 10      & 7       \\ \hline
2. & $c_r$    & 0.5      & 0.7     & 1.0     & 1.0     & 1.0      \\ \hline
3. & $c_T$    & 0.5      & 0.5     & 0.7     & 0.9     & 0.9       \\ \hline
4. & $q_r$    & 0.4      & 0.4     & 0.4     & 0.4     & 0.4        \\ \hline
5. & $q_T$    & 0.3      & 0.3     & 0.3     & 0.3     & 0.3         \\ \hline
6. & $\gamma$ & 0.1      & 0.1     & 0.1     & 0.15    & 0.21         \\ \hline
7. & $k_1$    & 0.00005  & 0.00005 & 0.00005 & 0.0005 & 0.0005       \\ \hline
8. & $S_{ra}$ & 6000     & 6000    & 6000    & 6000    & 6000    \\ 
\hline 
\end{tabular}}
\label{table1}
\end{table} 

\begin{figure}
    \centering 
        \includegraphics[width=0.6\textwidth]{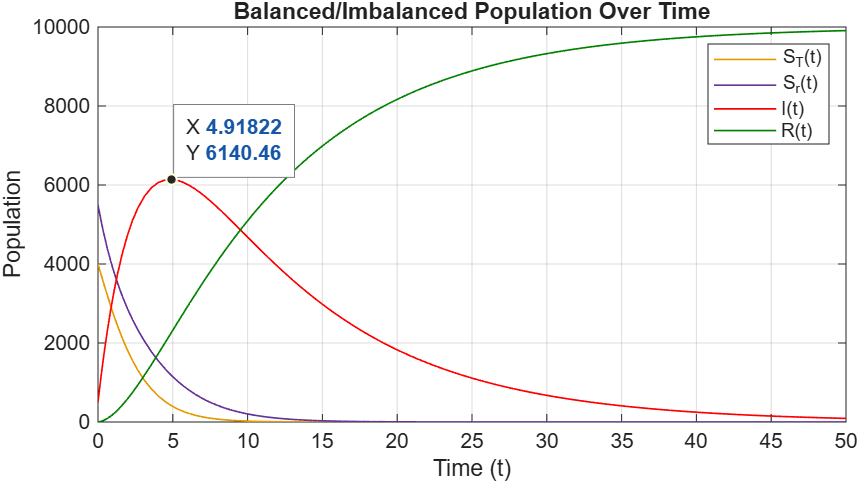}
        \caption{Simulation 1 Time evolution graph}
        \label{sim1}
\end{figure}

\begin{figure}[h!]
    \centering 
    \begin{subfigure}{0.49\linewidth}
        \centering
        \includegraphics[width=\textwidth]{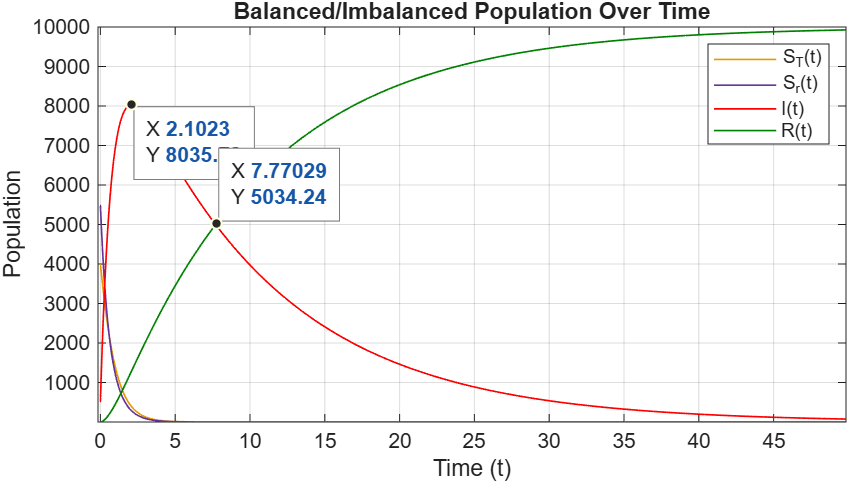}
        \subcaption{Simulation 2}
        \label{sim2}
    \end{subfigure}
    \hfill 
    \begin{subfigure}{0.49\linewidth}
        \centering
        \includegraphics[width=\textwidth]{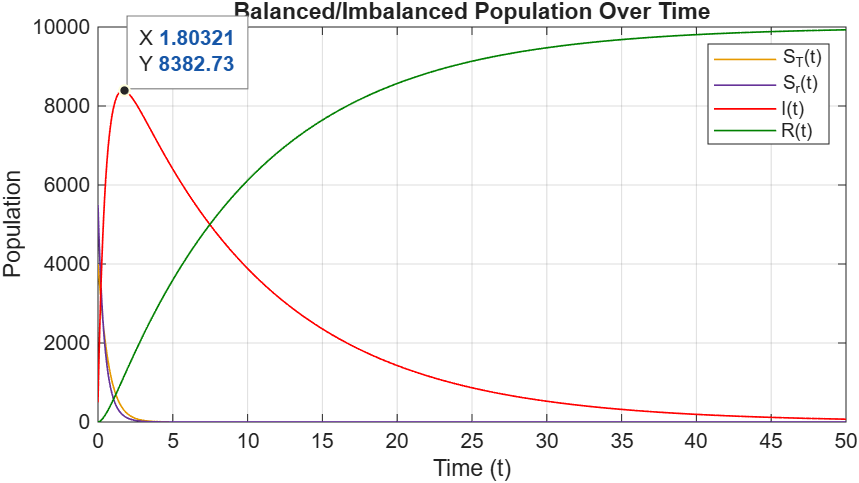}
        \subcaption{Simulation 3}
        \label{sim3}
    \end{subfigure}
    \caption{Time evolution graphs}
    \label{fig2}
\end{figure} 

\begin{figure}[h!]
    \centering 
    \begin{subfigure}{0.49\linewidth}
        \centering
        \includegraphics[width=\textwidth]{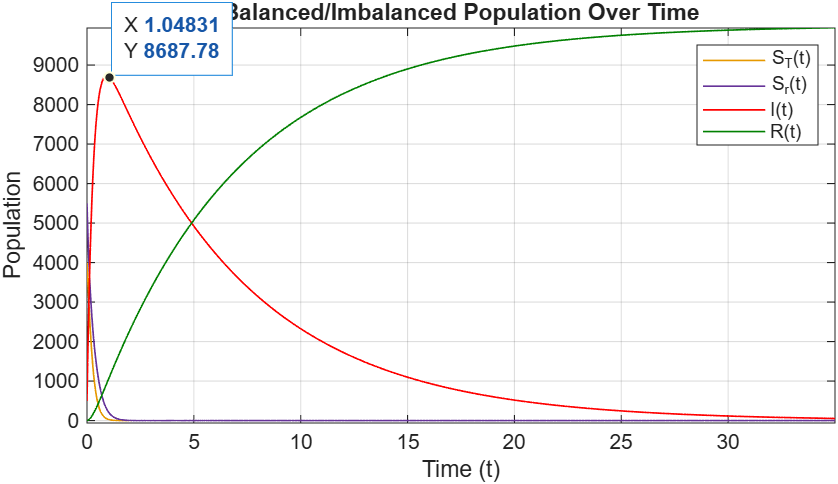}
        \subcaption{Simulation 4}
        \label{sim4}
    \end{subfigure}
    \hfill 
    \begin{subfigure}{0.49\linewidth}
        \centering
        \includegraphics[width=\textwidth]{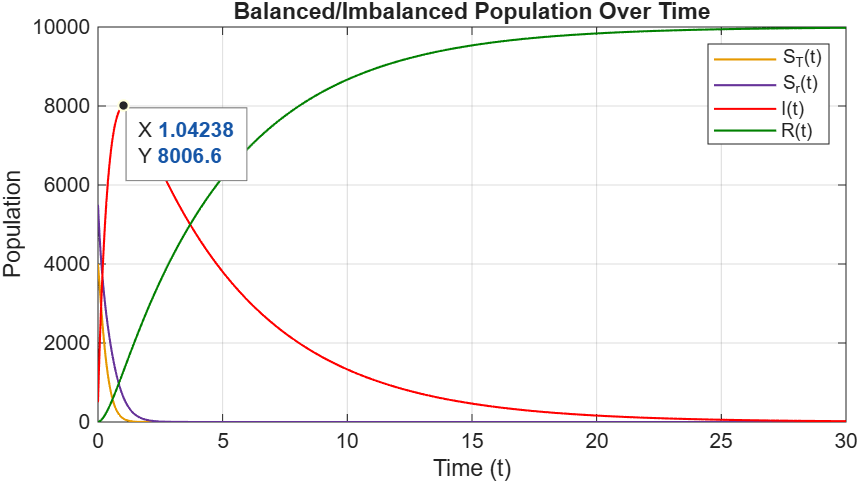}
        \subcaption{Simulation 5}
        \label{sim5}
    \end{subfigure}
    \caption{Time evolution graphs}
    \label{fig3}
\end{figure}

\begin{table}[h!]
\centering
\caption{Eigen Values corresponding to equilibrium $E_1$}
\resizebox{12 cm}{!}{
\begin{tabular}{|c|c|c|c|c|c|c|} 
\hline
\textbf{Simulation} & \textbf{Eigen Value 1 ($- nq_rc_r$)} & \textbf{Eigen Value 2 ($- nq_Tc_T$)} & \textbf{Eigen Value 3 ($- \gamma$)} \\ \hline
1. & -0.6 & -0.4 & -0.1            \\ \hline
2. & -1.2 & -1.68 & -0.1    \\ \hline
3.  & -1.56 & -2.4 & -0.1       \\ \hline
4. & -5.7 & -4.0 & -0.15        \\ \hline
5. & -4.8900 & -2.8 & -0.21         \\ \hline 
\end{tabular}}
\label{table3}
\end{table}

\begin{figure}[h!]
    \centering 
    \begin{subfigure}{0.49\linewidth}
        \centering
        \includegraphics[width=\textwidth]{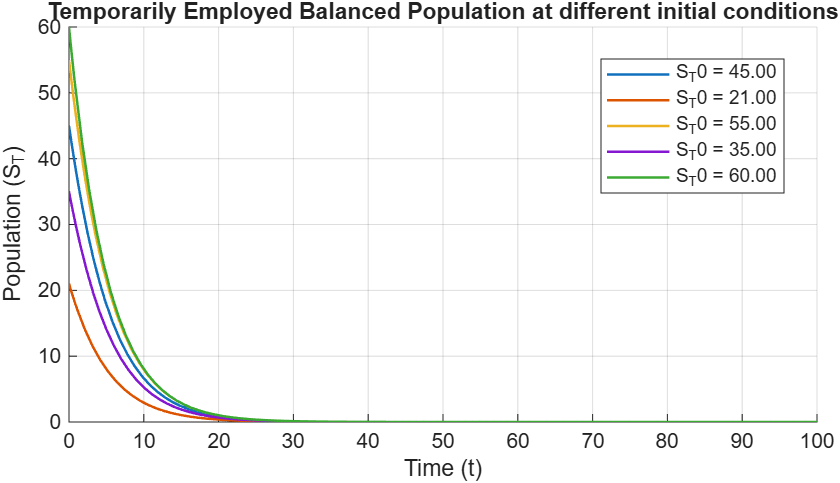}
        \subcaption{$S_T$}
        \label{s_T}
    \end{subfigure}
    \hfill 
    \begin{subfigure}{0.49\linewidth}
        \centering
        \includegraphics[width=\textwidth]{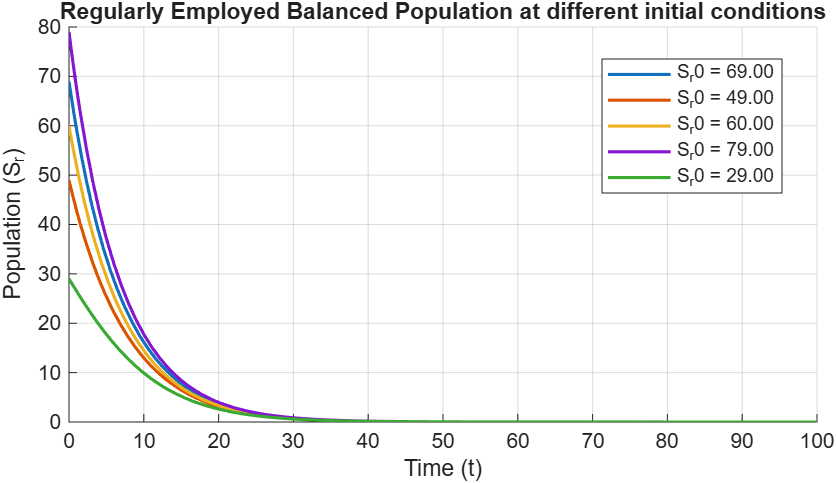}
        \subcaption{$S_r$}
        \label{s_r}
    \end{subfigure}
    \caption{Time evolution graphs with multiple initial conditions for $S_r$ and $S_T$}
    \label{fig4}
\end{figure}

\begin{figure}[h!]
    \centering 
    \begin{subfigure}{0.49\linewidth}
        \centering
        \includegraphics[width=\textwidth]{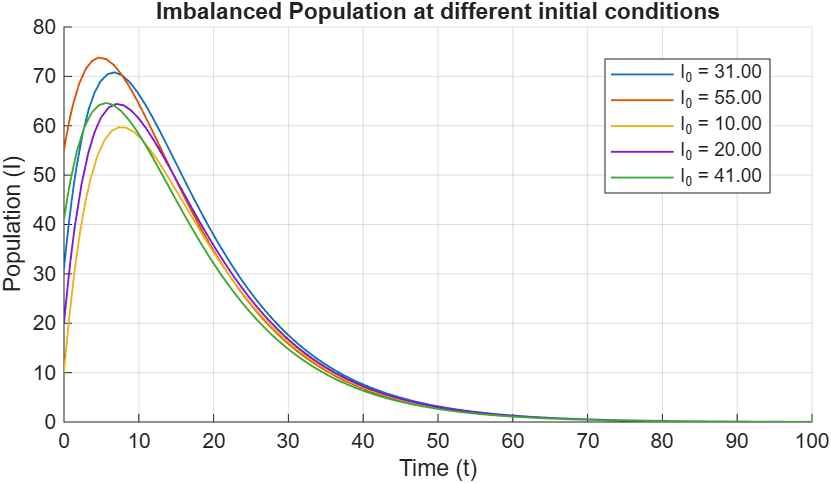}
        \subcaption{$I$}
        \label{I}
    \end{subfigure}
    \hfill 
    \begin{subfigure}{0.49\linewidth}
        \centering
        \includegraphics[width=\textwidth]{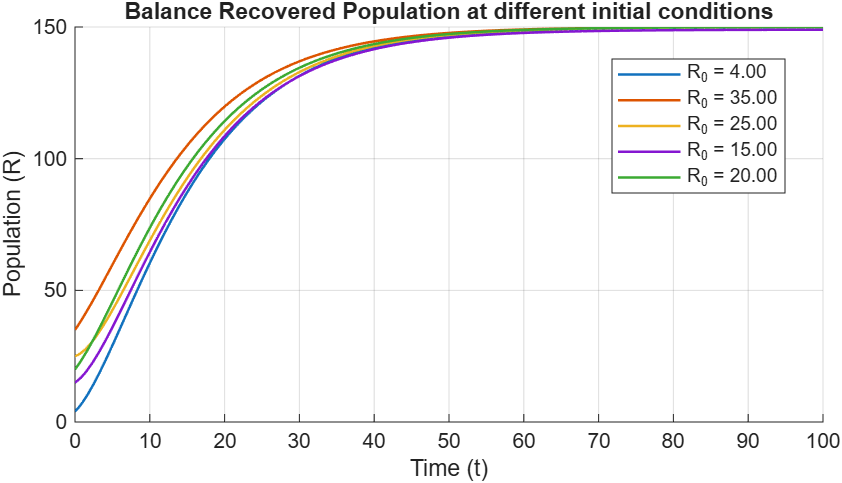}
        \subcaption{$R$}
        \label{R}
    \end{subfigure}
    \caption{Time evolution graphs with multiple initial conditions for $I$ and $R$}
    \label{fig5}
\end{figure}

\textit{Discussion: } It can be observed from the graphs that all the trajectories starting from different initial values settle at the equilibrium $E_1$ (0, 0, 0, N) over time for all simulations. Hence, it is inferred that inside the region of attraction, the WLB system is non-linearly stable around this equilibrium, that is, $E_1$. The eigenvalues corresponding to equilibrium $E_1$ are all negative (see from table \ref{table3}), making $E_1$ locally asymptotically stable. \newline 
On the other hand, the feasibility conditions for equilibrium $E_2$ are not satisfied by the set of parameter values in table \ref{table1} for simulations (1)-(5), making $E_2$ infeasible, meaning that $S_r^*$ and $R^*$ exceed the total population N. For simulation 1, it comes out to be (9000, 0, 13000), and using these values when the equilibrium $S_T^*$ for $E_2$ is calculated, it turns out to be negative. This violates the law of conservation of population in the closed system. One can compute $E_2$ for other simulations using Eq. (28) if required. 
\newline
For simulation 1, n = 2 has been taken, which represents the minimum number of roles on average a working individual must be engaged in, possibly one in non-work life towards self, and another as an employee in work life. \newline
From the figure \ref{sim1}, it is noticeable that the maximum population in the state of imbalance at one time point is just slightly over 6000. This number changes significantly and reaches over 8000 in \ref{sim2} with an increase in the value of $n$ and $c_r$. The curves have become steeper, while the ones for $S_r$ and $S_T$ almost overlapped each other. This showcases that for different initial values of $S_r$ and $S_T$, their trajectories become the same by increasing the values of $n$ and $c_r$. \newline

\begin{figure}[h!]
    \centering 
    \begin{subfigure}{0.49\linewidth}
        \centering
        \includegraphics[width=\textwidth]{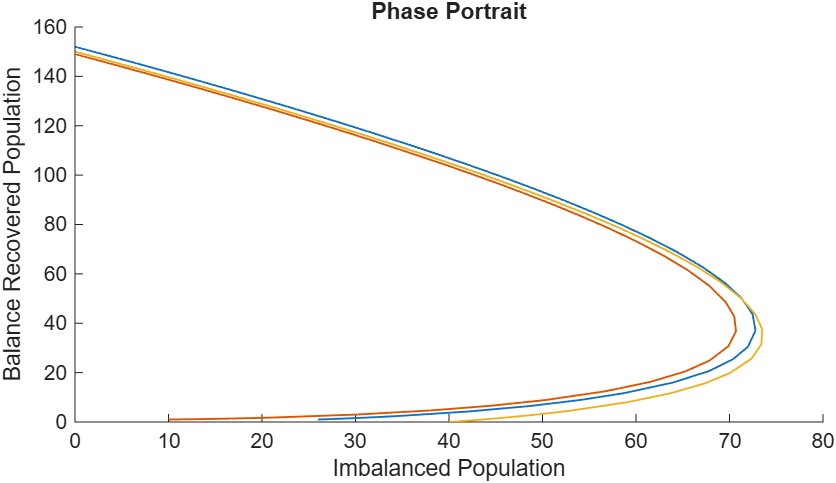}
        \subcaption{$I$ vs $R$}
        \label{pp5}
    \end{subfigure}
    \hfill 
    \begin{subfigure}{0.49\linewidth}
        \centering
        \includegraphics[width=\textwidth]{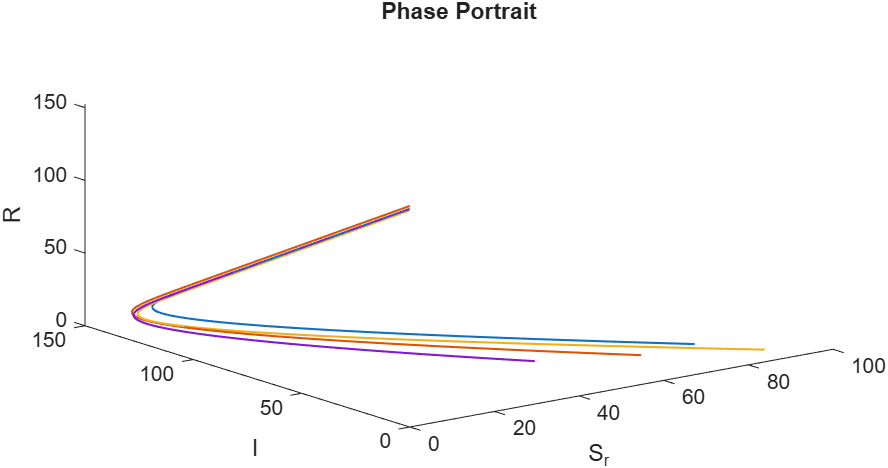}
        \subcaption{$S_r$ vs $I$ vs $R$}
        \label{pp6}
    \end{subfigure}
    \caption{Phase portrait of Imbalanced and Balance Recovered Individuals in 2D and 3D (with $S_r$)}
    \label{fig8}
\end{figure}

In all the simulations, one can observe that the dynamics become faster and faster as the parameter values are increased. This means that with an increase in the values of $n, c_r, c_T$, the population shifts faster to the imbalanced state, with the maximum number of individuals attaining the state at the same point in time. Furthermore, by increasing the value of $\gamma$ in simulations 4 and 5(\ref{sim4}, \ref{sim5}), the dynamics become even faster, and the system attains its equilibrium by time 30-35 units, whereas for other simulations it is post the 45 mark.

    \begin{figure}[h!]
    \centering 
    \begin{subfigure}{0.49\linewidth}
        \centering
        \includegraphics[width=\textwidth]{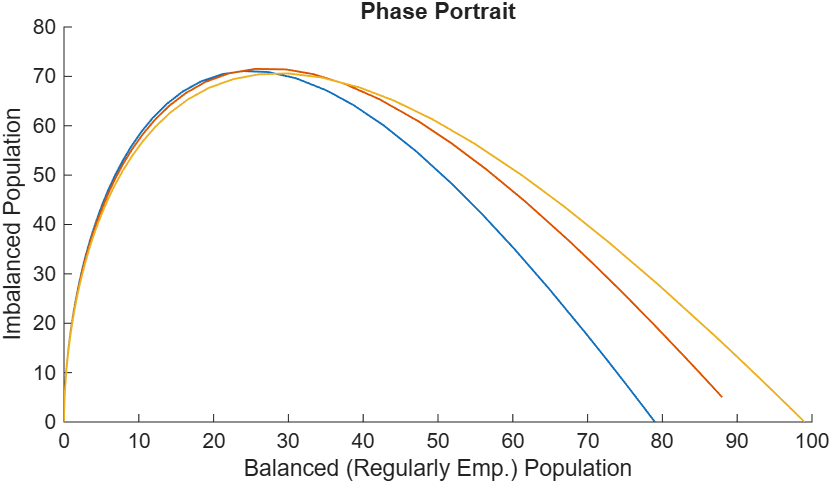}
        \subcaption{$S_r$ vs $I$}
        \label{pp1}
    \end{subfigure}
    \hfill 
    \begin{subfigure}{0.49\linewidth}
        \centering
        \includegraphics[width=\textwidth]{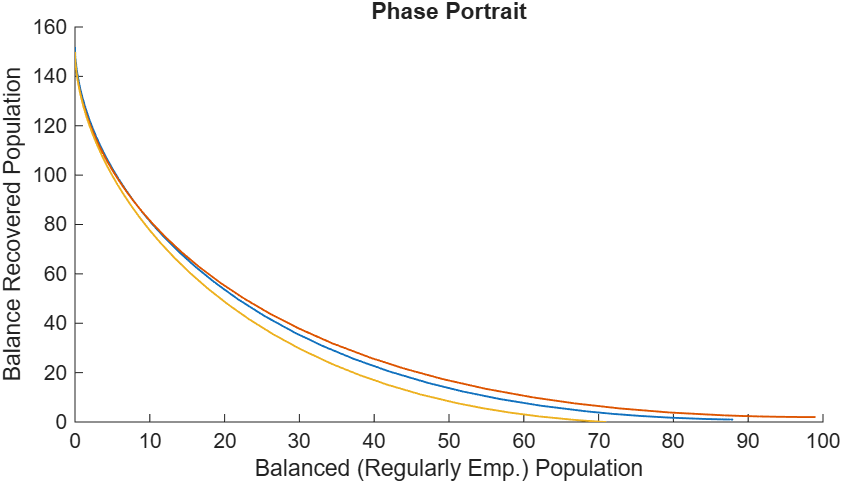}
        \subcaption{$S_r$ vs $R$}
        \label{pp2}
    \end{subfigure}
    \caption{Phase portraits of Regularly Employed Balanced Individuals}
    \label{fig6}
\end{figure}

\begin{figure}[h!]
    \centering 
    \begin{subfigure}{0.49\linewidth}
        \centering
        \includegraphics[width=\textwidth]{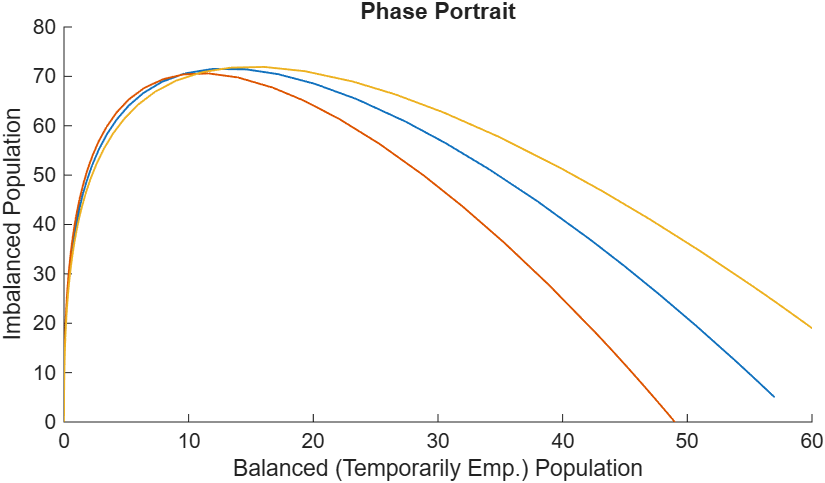}
        \subcaption{$S_T$ vs $I$}
        \label{pp3}
    \end{subfigure}
    \hfill 
    \begin{subfigure}{0.49\linewidth}
        \centering
        \includegraphics[width=\textwidth]{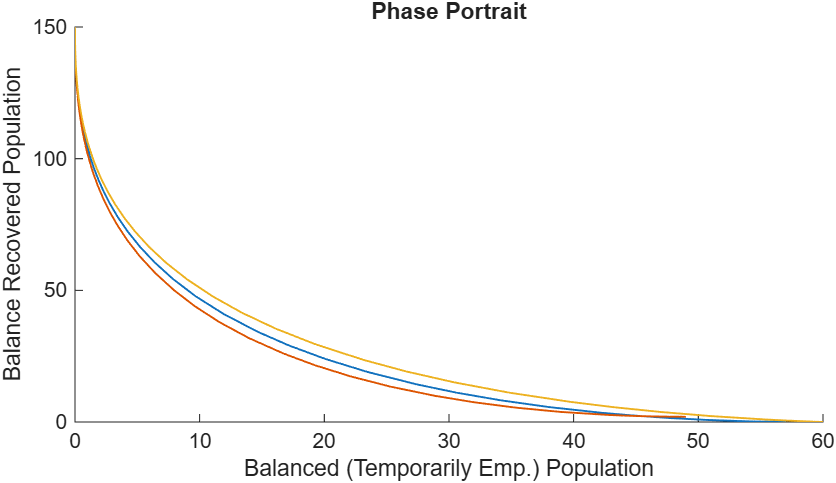}
        \subcaption{$S_T$ vs $R$}
        \label{pp4}
    \end{subfigure}
    \caption{Phase portraits of Temporarily Employed Balanced Individuals}
    \label{fig7}
\end{figure}

\section{Conclusion}
\label{Con} 

Through this study, we have established a foundational mathematical framework for understanding work/life balance as a dynamic, population-level phenomenon. By adapting compartmental modeling techniques from epidemiology, we have formalized how the working population transits between states of balance and imbalance. This interdisciplinary approach is achieved by integrating the sociological concepts of role theory and work/life conflict into a rigorous mathematical model. Using these concepts, we have created the measurable parameters that govern the state trajectories of the population. This reflects how our study bridges the existing gap between qualitative and quantitative analysis of WLB. \newline 

From our analysis, we have found two equillibrium points for the state of balance, out of which the only feasible and stable long-term outcome is a fully balanced population, represented by the disease-free equilibrium $E_1$. The outcome is achieved when recovery measures are adopted, and remain functional. On the other hand, the endemic equilibrium $E_2$, which is mathematically derivable, is observed to consistently violate population conservation under realistic parameterization. This indicates that the ongoing social imbalance may require to incorporate the model extensions. \newline

The results from numerical simulations offer actionable insights, that is, 

\begin{itemize}
    \item Reducing per-role complaints per unit time ($c_r, c_T$) directly increases the proportion of balanced individuals.
    \item Limiting the number of roles a member engages in ($n$) yields a similarly positive effect on the balanced population.
    \item Enhancing the recovery rate ($\gamma$) accelerates the dynamics of the transition and restoration of balance across the population.
\end{itemize}  

These findings quantify the intuitive notion that simplifying the role demands and strengthening coping resources are effective strategies for achieving ideal work/life balance for every individual. \newline

This research provides a robust starting point for future refinements, including demographic heterogeneity, stochastic effects, and empirical parameter estimation from survey data. The current model works under certain mathematical and theoretical assumptions. Primarily, the baseline model assumes demographic homogeneity across the population and posits that recovery from work/life imbalance is permanent. These assumptions, while necessary for establishing the initial framework, represent constraints on the model's direct applicability to highly complex, real-world populations. \newline 

As a future scope, the model can be refined by relaxing the current assumptions by incorporating demographic heterogeneity and stochastic effects to better reflect real-world population dynamics. There is a crucial need for empirical parameter estimation using real-world survey data to calibrate the model accurately. Moreover, the conflict parameter can further be bifurcated into three categories, discussed in chapter \ref{Pre}, to understand how much of WLB is affected by the amount of time devoted to work and non-work life, or other forms of conflict (strain and behavior-based), which also contribute significantly to the imbalance. Ultimately, this work bridges social theory and applied mathematics, offering a novel lens for diagnosing, predicting, and intervening in the growing crisis of work-life imbalance.

\section{Acknowledgements}

The authors thank the anonymous reviewers for thoughtful remarks and suggestions, which have helped to improve the manuscript significantly.

\bibliographystyle{IEEEtranN}
\bibliography{bibliography}

\begin{thebibliography}{13}
\providecommand{\natexlab}[1]{#1}
\providecommand{\url}[1]{#1}
\csname url@samestyle\endcsname
\providecommand{\newblock}{\relax}
\providecommand{\bibinfo}[2]{#2}
\providecommand{\BIBentrySTDinterwordspacing}{\spaceskip=0pt\relax}
\providecommand{\BIBentryALTinterwordstretchfactor}{4}
\providecommand{\BIBentryALTinterwordspacing}{\spaceskip=\fontdimen2\font plus
\BIBentryALTinterwordstretchfactor\fontdimen3\font minus \fontdimen4\font\relax}
\providecommand{\BIBforeignlanguage}[2]{{%
\expandafter\ifx\csname l@#1\endcsname\relax
\typeout{** WARNING: IEEEtranN.bst: No hyphenation pattern has been}%
\typeout{** loaded for the language `#1'. Using the pattern for}%
\typeout{** the default language instead.}%
\else
\language=\csname l@#1\endcsname
\fi
#2}}
\providecommand{\BIBdecl}{\relax}
\BIBdecl

\bibitem[Frone et~al.(1992)Frone, Russell, and Cooper]{frone1992antecedents}
M.~R. Frone, M.~Russell, and M.~L. Cooper, ``Antecedents and outcomes of work-family conflict: testing a model of the work-family interface.'' \emph{Journal of applied psychology}, vol.~77, no.~1, p.~65, 1992.

\bibitem[Meenakshi et~al.(2013)Meenakshi, Subrahmanyam, and Ravichandran]{meenakshi2013importance}
S.~P. Meenakshi, V.~Subrahmanyam, and K.~Ravichandran, ``The importance of work-life-balance,'' \emph{IOSR Journal of Business and Management}, vol.~14, no.~3, pp. 31--35, 2013.

\bibitem[Carlson et~al.(2000)Carlson, Kacmar, and Williams]{article}
D.~Carlson, K.~Kacmar, and L.~Williams, ``Construction and initial validation of a multidimensional measure of work–family conflict,,'' \emph{Journal of Vocational Behavior}, vol.~56, pp. 249--276, 04 2000.

\bibitem[G{\"u}lmez et~al.(2024)G{\"u}lmez, Koruca, Aydin, and Urganci]{gulmez2024heuristic}
E.~G{\"u}lmez, H.~I. Koruca, M.~E. Aydin, and K.~B. Urganci, ``Heuristic and swarm intelligence algorithms for work-life balance problem,'' \emph{Computers \& Industrial Engineering}, vol. 187, p. 109857, 2024.

\bibitem[Greenhaus and Beutell(1985)]{greenhaus1985sources}
J.~H. Greenhaus and N.~J. Beutell, ``Sources of conflict between work and family roles,'' \emph{Academy of management review}, vol.~10, no.~1, pp. 76--88, 1985.

\bibitem[Greenhaus and Powell(2006)]{greenhaus2006work}
J.~H. Greenhaus and G.~N. Powell, ``When work and family are allies: A theory of work-family enrichment,'' \emph{Academy of management review}, vol.~31, no.~1, pp. 72--92, 2006.

\bibitem[Fisher-McAuley et~al.(2003)Fisher-McAuley, Stanton, Jolton, and Gavin]{fisher2003modelling}
G.~Fisher-McAuley, J.~Stanton, J.~Jolton, and J.~Gavin, ``Modelling the relationship between work life balance and organisational outcomes,'' in \emph{Annual conference of the society for industrial-organisational psychology. Orlando}, vol.~1, no.~26, 2003, pp. 1--30.

\bibitem[Greenhaus and Powell(2003)]{GREENHAUS2003291}
J.~H. Greenhaus and G.~N. Powell, ``When work and family collide: Deciding between competing role demands,'' \emph{Organizational Behavior and Human Decision Processes}, vol.~90, no.~2, pp. 291--303, 2003.

\bibitem[McMillan et~al.(2011)McMillan, Morris, and Atchley]{mcmillan2011constructs}
H.~S. McMillan, M.~L. Morris, and E.~K. Atchley, ``Constructs of the work/life interface: A synthesis of the literature and introduction of the concept of work/life harmony,'' \emph{Human Resource Development Review}, vol.~10, no.~1, pp. 6--25, 2011.

\bibitem[Sirgy and Lee(2018)]{sirgy2018work}
M.~J. Sirgy and D.-J. Lee, ``Work-life balance: An integrative review,'' \emph{Applied Research in Quality of Life}, vol.~13, no.~1, pp. 229--254, 2018.

\bibitem[Misra and Singh(2011)]{misra2011mathematical}
A.~Misra and A.~K. Singh, ``A mathematical model for unemployment,'' \emph{Nonlinear Analysis: Real World Applications}, vol.~12, no.~1, pp. 128--136, 2011.

\bibitem[Vaidya et~al.(2026)Vaidya, Iqbal, Khan, Gupta, Ashirova, and Goyal]{vaidya2026analyzing}
G.~Vaidya, S.~M.~U. Iqbal, M.~I. Khan, K.~Gupta, A.~Ashirova, and S.~Goyal, ``Analyzing nonlinear dynamical systems and chaos using machine learning in physical and biological systems,'' \emph{Journal of Interdisciplinary Mathematics}, vol.~29, no.~3, pp. 663--670, 2026.

\bibitem[Kulkarni et~al.(2025)Kulkarni, Kadam, Paswan, Shipra, Julme, and Dhondge]{kulkarni2025numerical}
S.~Kulkarni, U.~Kadam, D.~K. Paswan, K.~Shipra, B.~C. Julme, and S.~L. Dhondge, ``Numerical stability analysis of nonlinear differential equations in engineering systems,'' \emph{Journal of Interdisciplinary Mathematics}, vol.~28, no.~6, pp. 2099--2108, 2025.

\end{thebibliography}

\end{document}